\documentclass[12pt]{article}
\usepackage{color}
\usepackage{amsrefs}
\usepackage{amsmath}
\usepackage{amsfonts}
\usepackage{amsthm}
\usepackage{geometry}
\usepackage{enumerate}
\ifdefined\directlua
\usepackage{polyglossia}
\setmainlanguage{english}

\usepackage{fontspec}
\usepackage{microtype}
\usepackage{pdftexcmds}
\makeatletter
\ifcase\pdf@shellescape
\relax\or
\begingroup
  \catcode`\%=12\relax
  \gdef\sformat{"
\endgroup
\directlua{
 local cmd="git show -s --format='"..\sformat.."'"
 print(cmd)
 local r=io.popen(cmd):read("*a")
 if (r) then
      tex.print("\string\\def\string\\COMMIT{"..r.."}")
 end
 }
\or
\relax\fi
\makeatother
\ifdefined\COMMIT
        \usepackage{background}
        \backgroundsetup{%
         pages=all, placement=bottom,angle=0,scale=2,%
         vshift=20pt,%
         contents={Commit version:\COMMIT}}
\fi

\fi

\newcommand{\wt}{\mathrm{wt}\,}
\newcommand{\Fix}{\mathrm{Fix}}
\newcommand{\FF}{\mathbb F}
\newcommand{\ZZ}{\mathbb Z}
\newcommand{\QQ}{\mathbb Q}
\newcommand{\CC}{\mathbb C}
\newcommand{\RR}{\mathbb R}
\newcommand{\SSS}{\mathbb S}
\newcommand{\KK}{\mathbb K}
\newcommand{\fq}{\mathfrak q}
\newcommand{\fb}{\mathfrak b}
\newcommand{\cD}{\mathcal D}
\newcommand{\cP}{\mathcal P}
\newcommand{\cV}{\mathcal V}
\newcommand{\cW}{\mathcal W}
\newcommand{\cS}{\mathcal S}
\newcommand{\cX}{\mathcal X}
\newcommand{\cY}{\mathcal Y}
\newcommand{\cZ}{\mathcal Z}
\newcommand{\cQ}{\mathcal Q}
\newcommand{\cL}{\mathcal L}
\newcommand{\cG}{\mathcal G}
\newcommand{\cT}{\mathcal T}
\newcommand{\GG}{\mathbb G}
\newcommand{\VV}{\mathbb V}
\newcommand{\LL}{\mathbb L}
\newcommand{\cC}{\mathcal C}
\newcommand{\cF}{\mathcal F}
\newcommand{\cH}{\mathcal H}
\newcommand{\Res}{\mathrm{Res}}
\newcommand{\GF}{\mathrm{GF}}
\newcommand{\Cp}{p}
\newcommand{\Cq}{q}
\newcommand{\Sec}{\mathrm{Sec}}
\newcommand{\PG}{\mathrm{PG}}
\newcommand{\PGL}{\mathrm{PGL}}
\newcommand{\GL}{\mathrm{GL}}
\newcommand{\fH}{\mathfrak H}
\newcommand{\fA}{\mathfrak A}
\newcommand{\fN}{\mathfrak N}
\newcommand{\fB}{\mathfrak B}
\newcommand{\fG}{\mathfrak G}
\newcommand{\fC}{\mathfrak C}
\newcommand{\cod}{\operatorname{codim}}
\newcommand{\diam}{\operatorname{diam}}
\newcommand{\Hom}{\mathrm{Hom}}
\newcommand{\Img}{\mathrm{Img}}
\newcommand{\Aut}{\mathrm{Aut}}
\newcommand{\rank}{\mathrm{rank}}
\newcommand{\rad}{\mathrm{Rad}}
\newtheorem{mainth}{Theorem}
\newtheorem{maincor}[mainth]{Corollary}
\newtheorem{mainprob}[mainth]{Problem}
\newtheorem{theorem}{Theorem}[section]
\newtheorem{lemma}[theorem]{Lemma}
\newtheorem{corollary}[theorem]{Corollary}
\newtheorem{prop}[theorem]{Proposition}
\theoremstyle{definition}

\newtheorem{note}{Note}

\newcommand{\fD}{\mathfrak D}
\newcommand{\fmm}{\mathfrak w}
\newcommand{\bS}{\mathbb S}
\title{Nearly all subspaces of a classical polar space arise from its universal embedding}
\author{I. Cardinali, L. Giuzzi and A. Pasini}
\begin{document}
\maketitle

\begin{abstract}
Let $\Gamma$ be an embeddable non-degenerate polar space of finite rank $n \geq 2$. Assuming that $\Gamma$ admits the universal embedding (which is true for all embeddable polar spaces except grids of order at least $5$ and certain generalized quadrangles defined over quaternion division rings), let $\varepsilon:\Gamma\to\PG(V)$ be the universal embedding of $\Gamma$. Let $\cal S$ be a subspace of $\Gamma$ and suppose that $\cal S$, regarded as a polar space, has non-degenerate rank at least $2$. We shall prove that $\cS$ is the $\varepsilon$-preimage of a projective subspace of $\PG(V)$.
\end{abstract}

\section{The main result of this paper}
\label{section1}
Let $\Gamma = (\cP, \cL)$ be a non-degenerate embeddable polar space of finite rank at least $2$, where $\cP$ and $\cL$ are the point-set and the line-set of $\Gamma$ respectively. It is well known that $\Gamma$ admits the (absolutely) universal embedding but for two exceptional cases of rank $n = 2$, where at least two distinct relatively universal embeddings exist (see Tits \cite[\S 8.6(II) (a), (b)]{Tits}, also cases (A) and (B) of Theorem \ref{Tits1} of this paper). Keeping those two exceptional cases aside, let $\varepsilon:\Gamma\to\PG(V)$ be the universal embedding of $\Gamma$.

We recall that a {\em subspace} of $\Gamma$ is a subset $\cS\subseteq\cP$ such that every line of $\Gamma$ meeting $\cS$ in two distinct points is fully contained in $\cS$. If $\cS \subset \cP$ then $\cS$ is called a {\em proper} subspace.

Let $\cS$ be a subspace of $\Gamma$. Then $\cS$, equipped with the lines of $\Gamma$ contained in it, is a possibly degenerate polar space. Let $\rad(\cS) = \cS^\perp\cap \cS$ be its radical and let $\rank(\rad(\cS))$ be the rank of $\rad(\cS)$, namely the size of a minimal spanning set of $\rad(\cS)$. We keep the symbol $\rank(\cS)$ for the rank of the polar space $\cS$ and put $\rank_{\mathrm{nd}}(\cS) := \rank(\cS)-\rank(\rad(\cS))$ ($= \rank(\cS)$ if $\cS$ is non-degenerate). Following Buekenhout and Cohen \cite{BC}, we call $\rank_{\mathrm{nd}}(\cS)$ the {\em non-degenerate rank} of $\cS$. Obviously, $\rank_{\mathrm{nd}}(\cS) = 0$ if and only if $\cS \subseteq \cS^\perp$, namely $\cS$ is a {\em singular subspace}.

We say that $\cS$ {\em arises from the embedding} $\varepsilon$ if $\varepsilon^{-1}(\langle\varepsilon(\cS)\rangle_V) = \cS$, where $\langle . \rangle_V$ stands for spans in $\PG(V)$ and, for a subset $X\subseteq \PG(V)$, we denote by $\varepsilon^{-1}(X)$ the $\varepsilon$-preimage  $\{x\in{\cP}~|~\varepsilon(x)\in X\}$ of $X$, as usual. Obviously, singular subspaces and the improper subspace arise from $\varepsilon$ (in fact they arise from any embedding of $\Gamma$). So, we are not going to consider them in the following theorem, which is our main result in this paper.

\begin{mainth}\label{main 1}
With $\Gamma$ and $\varepsilon$ as above, let $\cal S$ be a proper non-singular subspace of $\Gamma$ and suppose that $\rank_{\mathrm{nd}}(\cS) \geq 2$. Then $\cS$ arises from $\varepsilon$.
\end{mainth}

We shall prove Theorem \ref{main 1} in Sections~\ref{section3} and \ref{section4}.

\begin{note}
  The hypothesis $\rank_{\mathrm{nd}}(\cS) \geq 2$ cannot be removed from Theorem \ref{main 1}, as one can easily see by noticing that, when  $\rank_{\mathrm{nd}}(\cS) = 1$, then $\cS$ is just a collection of singular subspaces of rank $k+1$ containing a given singular subspace $R = \rad(\cS)$ of rank $k$ and such that no two of them are contained in a common singular subspace. In particular, if $\rank(\cS) = 1$ then $\cal S$ is a set of mutually non collinear points, which in general does not arise from any embedding. Even assuming that $\rank(\Gamma) = 2$ and $\cS$ is an ovoid of $\Gamma$, we cannot claim that $\cS$ arises from $\varepsilon$. Indeed several classical generalized quadrangles exist which admit non-classical ovoids.
\end{note}

\begin{note}\label{A and B}
The two embeddable cases which admit more than one relatively universal embedding have been excluded since the beginning of this section, but Theorem \ref{main 1} would be vacuous for them. Indeed neither of them admits proper non-singular subspaces of non-degenerate rank at least 2. This is obvious for one of them, namely  case (B) of Theorem \ref{Tits1}, since the generalized quadrangles considered in that case are grids, but the same is true in the other case too (case (A) of Theorem \ref{Tits1}), as we shall prove in Section~\ref{section5}.
\end{note}

Recall that a {\em hyperplane} of $\Gamma$ is a proper subspace $\cS$ such that every line of $\Gamma$ meets $\cS$ non-trivially. For instance, the perp $p^\perp$ of a point $p\in\cP$ is a hyperplane, called a {\em singular hyperplane}. It is well known that all hyperplanes of $\Gamma$ are maximal in the family of proper subspaces of $\Gamma$ (see e.g. Shult \cite[Lemma 7.5.1]{S}). Conversely, by exploiting Theorem \ref{main 1} we can prove the following:

\begin{maincor}\label{main 2}
Every maximal (proper) subspace of $\Gamma$ of rank at least $2$ is a hyperplane.
\end{maincor}
\noindent
{\bf Proof.} Let $\cal S$ be a maximal subspace of $\Gamma$ and suppose that $\rank(\cS) \geq 2$. Suppose firstly that $\cS$ is degenerate and let $p\in\rad(\cS)$. Then $\cS\subseteq p^\perp$. By maximality, $\cS = p^\perp$. On the other hand, let $\cS$ be non-degenerate. Then $\cS = \varepsilon^{-1}(\langle \varepsilon(\cS)\rangle_V)$ by Theorem \ref{main 1}. Accordingly, $\langle \varepsilon(\cS)\rangle_V$ is a proper subspace of $\PG(V)$. Let $H$ be a hyperplane of $\PG(V)$ containing $\langle \varepsilon(\cS)\rangle_V$. Then $\cH := \varepsilon^{-1}(H)$ is a hyperplane of $\Gamma$ and contains $\cS$. By maximality, $\cS = \cH$. \hfill $\Box$

\bigskip

Observe that a maximal singular subspace in general is neither a
hyperplane nor a singular hyperplane of $\Gamma$.

\begin{maincor}\label{main 3}
Let $\rank(\Gamma) = n > 2$. Then the hyperplanes of $\Gamma$ are precisely the maximal subspaces of $\Gamma$ of rank at least 2 (in fact, they have rank either $n-1$ or $n$).
\end{maincor}
\noindent
{\bf Proof.} This claim immediately follows from Corollary \ref{main 2} and the fact that all hyperplanes of $\Gamma$ have rank at least $n-1$.  \hfill $\Box$

\begin{note}
The following is known since long ago: if $\rank(\Gamma) > 2$ then all hyperplanes of $\Gamma$ arise from $\varepsilon$. A popular proof of this fact exploits a result of Ronan \cite[Section 1, Corollary 3]{Ron} and the fact that when $\rank(\Gamma) > 2$ all hyperplane complements of $\Gamma$ are simply connected (\cite[Section 3]{Pas90}, also Cohen and Shult \cite{CS90}). Now we can obtain it as a special case of Theorem \ref{main 1}, noticing that if $\rank(\Gamma) > 2$ then all hyperplanes of $\Gamma$ have non-degenerate rank at least 2.
\end{note}

\begin{mainprob}\label{main 4}
Suppose that $\rank(\Gamma) > 2$. Does $\Gamma$ admit maximal subspaces of rank $1$?
\end{mainprob}

\begin{mainprob}\label{main 5}
Let $\rank(\Gamma) = 2$. Is it true that every maximal subspace of $\Gamma$ of rank $1$ is a hyperplane (whence an ovoid)?
\end{mainprob}
Turning to generating sets, the following holds:

\begin{maincor}\label{main 6}
For $X\subseteq \cP$, suppose that $\rank_{\mathrm{nd}}(\langle X\rangle_\Gamma)\geq 2$. Then $X$ contains a subset $Y$ such that
$\langle Y\rangle_\Gamma = \langle X\rangle_\Gamma$ and $Y$ is minimal with respect to this property.
\end{maincor}
\noindent
{\bf Proof.} By Theorem \ref{main 1}, the embedding $\varepsilon$ induces a bijection between the family of subsets $Y\subseteq X$ such that $\rank_{\mathrm{nd}}(\langle Y\rangle_\Gamma) \geq 2$ and the family of subsets $Z\subseteq \varepsilon(X)$ such that the polar space $\langle Z\rangle_V\cap\varepsilon(\Gamma)$ has non-degenerate rank at least 2. The conclusion follows from well known properties of $\PG(V)$.  \hfill $\Box$

\begin{mainprob}\label{main 7}
Prove Corollary {\rm \ref{main 6}} in a synthetic way, without calling Theorem {\rm \ref{main 1}} for help.
\end{mainprob}
\begin{note}
The same conclusions as in Theorem \ref{main 1} and Corollary \ref{main 6} are obtained in \cite[Lemmas 2.3 and 2.6]{gr} under the additional hypothesis that the underlying division ring of $\Gamma$ is commutative and $\cS$ and $X$ contain a frame of $\Gamma$.
\end{note}

\paragraph{Organization of this paper.} Some basics on point-line geometries, projective embeddings and classical polar spaces are recalled in Section~\ref{section2}. Sections~\ref{section3} and~\ref{section4} contain the proof of Theorem \ref{main 1}. In the last section of this paper we prove that Theorem \ref{main 1} is indeed vacuous in the two exceptional cases mentioned in Note \ref{A and B}.

\section{A survey of basic notions and known facts}
\label{section2}
\subsection{Notation and terminology for projective spaces}

Given a vector space $V$, we use the symbol $\langle . \rangle_V$ to denote spans in $V$ as well as spans in the projective space $\PG(V)$ associated to $V$, thus avoiding the clumsy notation $\langle . \rangle_{\PG(V)}$. This notational ambiguity will be harmless.

Given a nonzero vector ${\bf v}\in V$, we put $[{\bf v}] := \langle {\bf v}\rangle_V$ and, given a subset $X\subseteq V$, we put $[X] := \{[{\bf x}]~|~ {\bf 0} \neq {\bf x}\in X\}$. Thus, $\langle [X]\rangle_V = [\langle X\rangle_V]$.

Given two vector spaces $V$ and $W$ over the same division ring, every semi-linear mapping $f:V\to W$ induces a mapping $[f]$ from  $\PG(V)\setminus[\ker(f)]$ to $\PG(W)$. Following Faure and Fr\"{o}licher \cite{FF} we call $[f]$ a {\em morphism} from $\PG(V)$ to $\PG(W)$ and $[\ker(f)]$ the {\em kernel} of $[f]$, also denoting it by the symbol $\ker([f])$. Clearly, a morphism $[f]$ is injective if and only if $\ker([f]) = \emptyset$. It is an {\em isomorphism} if it is both injective and surjective. The automorphisms of $\PG(V)$ are also called {\em collineations}.

For a subspace $X$ of $V$, we put $\PG(V)/[X] := \PG(V/X)$, calling it the {\em quotient} of $\PG(V)$ over $[X]$. If $p_X$ is the canonical projection of $V$ onto $V/X$, the projectivization $[p_X]$ of $p_X$ is called the {\em projection} of $\PG(V)$ onto $\PG(V)/[X]$.

\subsection{Notation and terminology for point-line geometries}

A point-line geometry is a pair $\Gamma = (\cP, \cL)$ where $\cP$ is a non-empty set and $\cL$ is a family of subsets of $\cP$, the elements of $\cP$ being called {\em points} and those of $\cL$ {\em lines}.

In Section~\ref{section1} we have freely used the symbol $\perp$, assuming that the reader knows its meaning. Anyway, we explain it here. When writing $p\perp q$ for two points $p$ and $q$ of $\Gamma$, we mean that $p$ and $q$ are {\em collinear}, namely they belong to a common line. We denote by $p^\perp$ the set of points collinear with a given point $p$, with $p\in p^\perp$ by convention. Given a subset $X\subseteq \cP$ we put $X^\perp := \cap_{p\in X}p^\perp$, also adopting obvious shortenings as $X^{\perp\perp}$ for $(X^\perp)^\perp$ or $\{X,Y\}^\perp$ for $(X\cup Y)^\perp$. The geometry $\Gamma$ is said to be {\em connected} if the graph $(\cP, \perp)$ is connected.

In Section~\ref{section1} we have explained what a subspace is. We have stated that definition for polar spaces, but it applies to arbitrary point-line geometries. We are not going to repeat it here. We only note that the intersection of a family of subspaces is still a subspace. In particular, the intersection of all subspaces of $\Gamma$ containing a given set of points $X\subseteq \cP$ is the smallest subspace containing $X$. We denote it by $\langle X \rangle_\Gamma$ and we call it the subspace {\em generated} (also {\em spanned}) by $X$. We can also describe $\langle X\rangle_\Gamma$ as follows: put $X_0 := X$ and, for every natural number $n$, put $X_{n+1} := X_n\cup\bigcup(\ell\in\cL~|~ |\ell\cap X_n| > 1)$. Then $\langle X\rangle_\Gamma = \cup_{n=0}^\infty X_n$.

Following Buekenhout and Cohen \cite{BC}, we say a subspace $\cS$ of $\Gamma$ is {\em singular} if $\cS \subseteq \cS^\perp$.

\subsection{Embeddings}

Let $\Gamma = (\cP, \cL)$ be a connected point-line geometry. We recall that a (full) {\em projective embedding} (henceforth called just {\em embedding}) of $\Gamma$ is an injective mapping $\varepsilon:\cP\to\PG(V)$ for some vector space $V$, such that $\varepsilon(\cP)$ spans $\PG(V)$ and the set $\varepsilon(\ell) = \{\varepsilon(x)\}_{x\in\ell}$ is a line of $\PG(V)$ for every line $\ell\in{\cL}$. We denote by $\varepsilon(\Gamma)$ the $\varepsilon$-image of $\Gamma$, namely the subgeometry of $\PG(V)$ with $\varepsilon(\cP)$ as the point-set and the projective lines $\varepsilon(\ell)$ as lines, for $\ell\in\cL$.

If $\KK$ is the underlying division of $V$ then we say that $\varepsilon$ is {\em defined over} $\KK$. If all embeddings of $\Gamma$ are defined over the same division ring $\KK$ then we say that $\Gamma$ is {\em defined over} $\KK$, also that $\KK$ is the {\em underlying division ring} of $\Gamma$.

We say that $\Gamma$ is {\em embeddable} if it admits at least one embedding. Note that $\Gamma$ is embeddable only if no two distinct lines of $\Gamma$ have two distinct points in common, as it follows from the fact that embeddings are injective and map lines of $\Gamma$ surjectively onto projective lines.

\subsubsection{Morphisms of embeddings}

Given two embeddings $\varepsilon_1:\Gamma\to\PG(V_1)$ and $\varepsilon_2:\Gamma\to\PG(V_2)$ defined over the same division ring, a {\em morphism} ({\em isomorphism}) from $\varepsilon_1$ to $\varepsilon_2$ is a morphism (isomorphism) $\varphi:\PG(V_1)\rightarrow\PG(V_2)$ such that $\varepsilon_2 = \varphi\cdot\varepsilon_1$. The subspace $K := \ker(\varphi)$ is such that $K\cap \varepsilon_1(\cP) = \emptyset$ and every line of $\PG(V_1)$ joining two distinct points of $\varepsilon_1(\cP)$ meets $K$ trivially. Conversely, for every subspace $K$ of $\PG(V_1)$ satisfying these properties, if $\pi_K$ is the projection of $\PG(V_1)$ onto $\PG(V_1)/K$ then the composite $\varepsilon_1/K := \pi_K\cdot\varepsilon_1$ is an embedding of $\Gamma$ and $\pi_K$ is a morphism from $\varepsilon_1$ to it. We call $\varepsilon_1/K$ the {\em quotient} of $\varepsilon_1$ over $K$. Clearly, if $\varphi:\varepsilon_1\to\varepsilon_2$ is a morphism of embeddings then $\varepsilon_2$ and $\varepsilon_1/\ker(\varphi)$ are isomorphic. By a little abuse, we say that $\varepsilon_2$ is a {\em quotient} of $\varepsilon_1$ over $\ker(\varphi)$.

The connectedness of $\Gamma$ implies that a morphism $\varphi:\varepsilon_1\to\varepsilon_2$, if it is exists, is unique (see e.g. Pasini and Van Maldeghem \cite[Section 4.1]{PVM}). Accordingly, if a morphism exists from $\varepsilon_1$ to $\varepsilon_2$ we write $\varepsilon_1 \to \varepsilon_2$, with no explicit mention of that morphism. We denote isomorphism of embeddings by the standard symbol $\cong$. Note that the uniqueness of the morphism between two embeddings implies that the identity mapping of $\PG(V)$ is the unique morphism from an embedding $\varepsilon:\Gamma\to\PG(V)$ to itself. Accordingly, if $\varepsilon_1\to\varepsilon_2\to\varepsilon_1$ for two embeddings $\varepsilon_1$ and $\varepsilon_2$ of $\Gamma$, then $\varepsilon_1\cong\varepsilon_2$.

\subsubsection{Universality}

An embedding $\varepsilon$ of $\Gamma$ is said to be {\em relatively universal} if $\varepsilon' \cong \varepsilon$ for any embedding $\varepsilon'$ of $\Gamma$ such that $\varepsilon'\to \varepsilon$. If moreover $\varepsilon \to \varepsilon'$ for any embedding $\varepsilon'$ of $\Gamma$ then $\varepsilon$ is said to be {\em absolutely universal} (also just {\em universal}, for short). Clearly, the absolutely universal embedding, if it exists, is unique up to isomorphisms.

Not every embeddable geometry admits the absolutely universal embedding. If a geometry admits it then it also admits an underlying division ring, but the latter property is not sufficient for the absolutely universal embedding to exist. On the other hand, every embeddable geometry admits relatively universal embeddings. Indeed, for every embedding $\varepsilon$ of $\Gamma$ there exists an embedding $\tilde{\varepsilon}\to\varepsilon$, sometimes called the {\em hull} of $\varepsilon$, such that $\tilde{\varepsilon}\to\varepsilon'$ for every embedding $\varepsilon'\to\varepsilon$ (Ronan \cite{Ron}). Clearly, $\tilde{\varepsilon}$ is uniquely determined by $\varepsilon$ modulo isomorphisms and it is relatively universal. Accordingly, the embeddings of $\Gamma$ are partitioned in mutually disjoint families where each family consists of all quotients of a given relatively universal embedding. The geometry $\Gamma$ admits the absolutely universal embedding if and only if just one such family exists, namely $\Gamma$ admits just one relatively universal embedding (modulo isomorphisms, of course).

\subsubsection{Homogeneity}\label{lift}

Given an embedding $\varepsilon:\Gamma\to\PG(V)$, we say that an automorphism $g\in \mathrm{Aut}(\Gamma)$ {\em lifts} to $\PG(V)$ {\em through} $\varepsilon$ if there exists a collineation $\varepsilon(g)$ of $\PG(V)$ such that the composite $\varepsilon\cdot g$ is an embedding of $\Gamma$ such that $\varepsilon\cdot g = \varepsilon(g)_{|\varepsilon(\cP)}\cdot\varepsilon$. Clearly $\varepsilon(g)$ stabilizes $\varepsilon(\cP)$ and induces on $\varepsilon(\cP)$ an automorphism of $\varepsilon(\Gamma)$. In short, it {\em stabilizes} $\varepsilon(\Gamma)$. The connectedness of $\Gamma$ implies that $\varepsilon(g)$, if it exists, is uniquely determined by $g$. We call $\varepsilon(g)$ {\em the lifting} of $g$ to $\PG(V)$. We also denote by $\mathrm{Aut}_\varepsilon(\Gamma)$ the subgroup of $\mathrm{Aut}(\Gamma)$ formed by the automorphisms of $\Gamma$ which lift to $\PG(V)$ via $\varepsilon$. Clearly, $\varepsilon(\mathrm{Aut}_\varepsilon(\Gamma)) = \{\varepsilon(g)~|~g\in \mathrm{Aut}_\varepsilon(\Gamma)\}$ is the stabilizer of $\varepsilon(\Gamma)$ in the collineation group of $\PG(V)$. The embedding $\varepsilon$ is said to be {\em homogeneous} if $\mathrm{Aut}_\varepsilon(\Gamma) = \mathrm{Aut}(\Gamma)$.

We have $\varepsilon\cdot g \cong \varepsilon$ if and only if $g\in \mathrm{Aut}_\varepsilon(\Gamma)$. It follows that all absolutely universal embeddings are homogeneous.

\subsubsection{Embeddings and subspaces}

Given a geometry $\Gamma = (\cP, \cL)$ and an embedding $\varepsilon:\Gamma\to\PG(V)$, the $\varepsilon$-preimage $\varepsilon^{-1}(\cX)$ of a projective subspace $\cX$ of $\PG(V)$ is a subspace of $\Gamma$. In particular, $\varepsilon^{-1}(\langle\varepsilon(X)\rangle_V)$ is a subspace of $\Gamma$ for every set of points $X\subseteq \cP$. Since it contains $\varepsilon^{-1}(\varepsilon(X)) = X$, we obtain that $\langle X\rangle_\Gamma  \subseteq   \varepsilon^{-1}(\langle \varepsilon(X)\rangle_V).$

As recalled in Section~\ref{section1}, a subspace $\cS$ of $\Gamma$ is said to {\em arise from} $\varepsilon$ if it is the $\varepsilon$-preimage of a subspace of $\PG(V)$; equivalently, $\varepsilon(\cS) = \langle\varepsilon(\cS)\rangle_V\cap\varepsilon(\cP)$. Clearly, if $\cS$ arises from $\varepsilon$ and $\varepsilon'\to\varepsilon$ then $\cS$ also arises from $\varepsilon'$. So, if $\Gamma$ admits the absolutely universal embedding then all subspaces of $\Gamma$ which arise from some embedding also arise from the universal one.

\subsection{Embeddable polar spaces}

We are not going to recall all basics on polar spaces here. We refer to Buekenhout and Cohen \cite[Chapters 7, 8]{BC} for this matter. In this section we shall only describe the embeddings of the classical (namely embeddable) polar spaces. Our main sources for this topic are Tits \cite[Chapter 8]{Tits} and Buekenhout and Cohen \cite[Chapters 9, 10]{BC}.

\subsubsection{Reflexive sesquilinear forms}\label{rsf}

Given a division ring $\KK$, an {\em admissible pair} for $\KK$ is a pair $(\sigma,\epsilon)$ where $\sigma$ is an anti-automorphism of $\KK$, $\epsilon\in \KK^* := \KK\setminus\{0\}$, $\epsilon^\sigma = \epsilon^{-1}$ and $t^{\sigma^2} = \epsilon t\epsilon^{-1}$ for every $t\in \KK$. Note that $\sigma^2 = \mathrm{id}_\KK$ if and only if $\epsilon$ belongs to the center of $\KK$.

Given a $\KK$-vector space $V$ and an admissible pair $(\sigma,\epsilon)$ of $\KK$, a $(\sigma,\epsilon)$-{\em sesquilinear form} is a mapping $f:V\times V\to \KK$ such that $f({\bf z}, {\bf x}s+{\bf y}t) =  f({\bf z}, {\bf x})s + f({\bf z},{\bf y})t$ and $f({\bf y}, {\bf x}) = f({\bf x},{\bf y})^\sigma\epsilon$ for all ${\bf x}, {\bf y}, {\bf z}\in V$ and $s, t\in \KK$ (consequently $f({\bf x}s+{\bf y}t,{\bf z}) = s^\sigma f({\bf x},{\bf z}) + t^\sigma f({\bf y},{\bf z})$).  A {\em reflexive sesquilinear form} is a $(\sigma,\epsilon)$-sesquilinear form, for some admissible pair $(\sigma,\epsilon)$.

Let $f$ be a $(\sigma,\epsilon)$-sesquilinear form with $\sigma \neq \mathrm{id}_\KK$ and $\epsilon\in\{1,-1\}$. If $\epsilon  = 1$ ($\epsilon = -1$) then $f$ is called {\em hermitian} ({\em anti-hermitian}). On the other hand, let $\sigma = \mathrm{id}_\KK$. Then $\KK$ is commutative and necessarily $\epsilon \in \{1, -1\}$. If $\epsilon = 1$ the form $f$ is called {\em symmetric}. If $f({\bf x},{\bf x}) = 0$ for every ${\bf x}\in V$ (as it is the case when $\mathrm{char}(\KK)\neq 2$ and $(\sigma,\epsilon) = (\mathrm{id}_\KK, -1)$) then $f$ is called {\em alternating}.

The {\em radical} $\rad(f)$ of a reflexive sesquilinear form $f$ is the set of vectors ${\bf v}\in V$ such that $f({\bf v},{\bf x}) = 0$ for every ${\bf x}\in V$. The form $f$ is said to be {\em non-degenerate} if $\rad(f) = \{{\bf 0}\}$.

A subspace $X\subseteq V$ is said to be {\em totally~isotropic} for $f$ (also {\em totally $f$-isotropic} for short) if $f({\bf x},{\bf y}) = 0$ for any  choice of ${\bf x},{\bf y}\in X$. Similarly, a vector ${\bf v}\in V$ is said to be {\em isotropic} for $f$ (also {\em $f$-isotropic}) if $f({\bf v},{\bf v}) = 0$. Accordingly, a subspace $[X]$ (a point $[{\bf v}]$) of $\PG(V)$ is {\em totally $f$-isotropic} ($f$-{\em isotropic}) if $X$ is totally $f$-isotropic ($\bf v$ is isotropic).

The form $f$ is {\em trace-valued} if the $f$-isotropic vectors span $V$. This is always the case when $\mathrm{char}(\KK) \neq 2$ (Tits \cite[8.1.6]{Tits}).

The $f$-isotropic points and the totally $f$-isotropic lines of $\PG(V)$ form a polar space $\Gamma(f)$, the singular subspaces of which are the totally $f$-isotropic subspaces of $\PG(V)$. In particular, the projective radical $[\rad(f)]$ of $f$ is the radical of $\Gamma(f)$. So, the polar space $\Gamma(f)$ is non-degenerate if and only if $f$ is non-degenerate. Assuming that $f$ is trace-valued, the natural inclusion mapping yields and embedding of $\Gamma(f)$ in $\PG(V)$.

Let $f$ be a $(\sigma,\epsilon)$-sesquilinear form and, given $\kappa\in \KK^*$, let $g := \kappa f$. Then $g$ is a $(\rho,\eta)$-sesquilinear form where $t^\rho = \kappa t^\sigma \kappa^{-1}$ for every $t\in \KK$ and $\eta = \kappa\kappa^{-\sigma}\epsilon$. We say that $g$ and $f$ are {\em proportional}. If $\sigma\neq \mathrm{id}_\KK$ we can always choose $\kappa$ in such a way that $g$ is hermitian or anti-hermitian, as we prefer \cite[8.1.2]{Tits}. In other words, a reflexive sesquilinear form which is neither symmetric nor alternating is always proportional to a suitable hermitian form as well as an anti-hermitian form.

Clearly, if $f$ and $g$ are proportional then $\Gamma(f) = \Gamma(g)$. Conversely, if $\Gamma(f) = \Gamma(g)$ then $f$ and $g$ are proportional. This is not so difficult to prove when $f$ and $g$ are bilinear or hermitian; of course, if one of them is bilinear but the other one is hermitian then  $\Gamma(f) \neq \Gamma(g)$. When either $f$ or $g$ is neither bilinear nor hermitian, we can replace it with a hermitian form proportional to it.

\subsubsection{Pseudoquadratic forms}\label{pseudo}

With $(\sigma,\epsilon)$ as above, put $\KK_{\sigma,\epsilon} := \{t-t^\sigma\epsilon\}_{t\in\KK}$. This is a subgroup of the additive group of $\KK$.  Regarded $\KK$ as a group, let $\overline{\KK}_{\sigma,\epsilon} := \KK/\KK_{\sigma,\epsilon}$. Note that $t^\sigma\KK_{\sigma,\epsilon}t\subseteq \KK_{\sigma,\epsilon}$ for every $t\in \KK$. So, for $t, s\in \KK$ and $\bar{s} = s+\KK_{\sigma,\epsilon} \in \overline{\KK}_{\sigma,\epsilon}$, the symbol $t^\sigma\bar{s}t$ denotes a well defined element of $\overline{\KK}_{\sigma,\epsilon}$.

A $(\sigma,\epsilon)$-quadratic form is a mapping $Q:V\to\overline{\KK}_{\sigma,\epsilon}$ such that
\begin{equation}\label{quad1}
Q({\bf x}t) ~ = ~ t^\sigma Q({\bf x})t, \hspace{5 mm} (\forall {\bf x}\in V, t\in \KK)
\end{equation}
and there exists a trace-valued $(\sigma,\epsilon)$-sesquilinear form $f_Q:V\times V\rightarrow\KK$ such that
\begin{equation}\label{quad2}
Q({\bf x}+{\bf y}) ~ = ~ Q({\bf x})+ Q({\bf y}) + (f_Q({\bf x},{\bf y})+\KK_{\sigma,\epsilon}), \hspace{5 mm} (\forall {\bf x}, {\bf y}\in V).
\end{equation}
A {\em pseudoqudratic form} is a $(\sigma,\epsilon)$-quadratic form, for some admissible pair $(\sigma,\epsilon)$. When $\KK_{\sigma, \epsilon} = \{0\}$ (equivalently, $(\sigma, \epsilon)= (\mathrm{id}_\KK, 1)$) we get back the usual {\em quadratic forms}.

The sesquilinear form $f_Q$ satisfying (\ref{quad2}) is uniquely determined by $Q$ except when the group $\overline{\KK}_{\sigma,\epsilon}$ is trivial, namely $\KK_{\sigma.\epsilon} = \KK$ (Tits \cite[8.2.3]{Tits}). This happens precisely when $(\sigma,\epsilon) = (\mathrm{id}_{\KK}, -1)$ and $\mathrm{char}(\KK) \neq 2$, but we exclude this case. So, for the rest of this subsection we assume that either $(\sigma, \epsilon) \neq (\mathrm{id}_{\KK}, -1)$ or $\mathrm{char}(\KK) = 2$. Accordingly, $f_Q$ is uniquely determined by $Q$ and $Q$ is non-trivial. We call $f_Q$ the {\em sequilinearization} of $Q$ (also the {\em bilinearization} of $Q$ when $Q$ is quadratic).

The {\em radical} $\rad(Q)$ of $Q$ is the set $\rad(Q) := Q^{-1}(\bar{0})\cap\rad(f_Q)$, where $\bar{0} := \KK_{\sigma,\epsilon}$ is the null element of the group $\overline{\KK}_{\sigma,\epsilon}$. This set is a subspace of $V$. The form $Q$ is said to be {\em non-degenerate} (also {\em non-singular}) if $\rad(Q) = \{{\bf 0}\}$.

A vector ${\bf v}\in V$ is said to be {\em singular} for $Q$ (also $Q$-{\em singular}) if $Q({\bf x}) = \bar{0}$ and a subspace $X\subseteq V$ is {\em totally singular} for $Q$ (also {\em totally $Q$-singular}) if all of its vectors are $Q$-singular. The same terminology is used for points and subspaces of $\PG(V)$.

The $Q$-singular points and the totally $Q$-singular lines of $\PG(V)$ form a polar space $\Gamma(Q)$, the singular subspaces of which are the totally $Q$-singular subspaces of $\PG(V)$. The projective radical $[\rad(Q)]$ of $Q$ is the radical of $\Gamma(Q)$. So, $\Gamma(Q)$ is non-degenerate if and only if $Q$ is non-degenerate. All $Q$-singular points of $\PG(V)$ are $f_Q$-isotropic and all totally $Q$-singular subspaces of $\PG(V)$ are totally $f_Q$-isotropic. Thus, $\Gamma(Q)$ is a subgeometry of $\Gamma(f_Q)$ (but not a subspace, in general). In particular, $[\rad(Q)]\subseteq [\rad(f_Q)]$.

When $\mathrm{char}(\KK)\neq 2$ the pseudoquadratic form $Q$ is uniquely determined by its sesquilinearization $f_Q$ and we have $\Gamma(Q) = \Gamma(f_Q)$ (Tits \cite[8.2.4]{Tits}). In this case we can safely ignore pseudoquadratic forms, if we like. On the other hand, when
$\mathrm{char}(\KK) = 2$, in general the polar space $\Gamma(Q)$ is a proper subgeometry of $\Gamma(f_Q)$. In this case it can happen that  $[\rad(Q)]\subseteq [\rad(f_Q)]$, as when $Q$ is non-degenerate but $f_Q$ is degenerate.

Proportionality can be defined for pseudoquadratic forms too. Let $\kappa\in \KK^*$. It is not difficult to check that $\kappa \KK_{\sigma, \epsilon} = \KK_{\rho,\eta}$, with $t^\rho = \kappa t^\sigma \kappa^{-1}$ and $\eta = \kappa\kappa^{-\sigma}\epsilon$ as in the previous subsection. Accordingly, left multiplication by $\kappa$ yields  a group isomorphism $\iota_\kappa: \overline{\KK}_{\sigma,\epsilon}\to\overline{\KK}_{\rho, \eta}$. The composite $\iota_\kappa\cdot Q$ is a $(\rho,\eta)$-quadratic form. We denote it by $\kappa Q$ for short and we say that $Q$ and $\kappa Q$ are {\em proportional}. It is clear that $\Gamma(\kappa Q) = \Gamma(Q)$, but the converse also holds true: if $\Gamma(Q') = \Gamma(Q)$ for two pseudoquadratic forms $Q$ and $Q'$, then $Q$ and $Q'$ are proportional.

Likewise in the case of reflexive sesquilinear forms, if $\sigma\neq \mathrm{id}_\KK$ we can always choose $\kappa$ in such a way that $\eta$ is either $1$ or $-1$, as we like.

\subsubsection{Generalized pseudoquadratic forms}\label{gp}

Reflexive sesquilinear forms and pseudoquadratic forms are enough to describe (relatively) universal embeddings of polar spaces, but not all embeddings of polar spaces are universal. In order to describe non-universal embeddings we need a generalization of pseudoquadratic forms.

We firstly discuss some properties of the groups $\KK_{\sigma,\epsilon}$ and $\overline{\KK}_{\sigma,\epsilon}$. For every element $t$ of the additive group $\KK$ and every element $s$ of the division ring $\KK$, we put $t\circ s := s^\sigma t s$. The `multiplication' $\circ$ defined in this way enjoys the usual properties of multiplication of vectors and scalars except right distributivity, which holds if and only if $\sigma = \mathrm{id}_\KK$ and $\mathrm{char}(\KK) = 2$.

A subgroup $R$ of $\KK$ is $\circ$-{\em closed} if $R\circ\KK \subseteq R$. For instance, both $\KK_{\sigma,\epsilon}$ and the group $\KK^{\sigma,\epsilon} := \{t\in\KK~|~ t+t^\sigma\epsilon = 0\} \supseteq \KK_{\sigma,\epsilon}$ are $\circ$-closed. If $R$ is $\circ$-closed and contains $\KK_{\sigma,\epsilon}$ then $\circ$ defines a `vector-by-scalar' multiplication in the quotient group $\overline{\KK}_R := \KK/R$, which we still denote by the symbol $\circ$. In particular, a vector-by-scalar multiplication $\circ$ is thus defined in $\overline{\KK}_{\sigma,\epsilon}$.
The right distributive property holds for an element $\bar{t} = t + \KK_{\sigma,\epsilon} \in \overline{\KK}_{\sigma,\epsilon}$ with respect to $\circ$ if and only if $t\in \KK^{\sigma,\epsilon}$. Thus, the quotient group $\overline{\KK}_{\sigma,\epsilon}^\circ := \KK^{\sigma,\epsilon}/\KK_{\sigma,\epsilon}$ is the largest subgroup $\overline{R}$ of $\overline{\KK}_{\sigma,\epsilon}$ such that $(\overline{R},\circ)$ is a $\KK$-vector space (see e.g. \cite[2.1.2]{Pas}).

Let $R$ be a $\circ$-closed proper subgroup of $\KK$ containing $\KK_{\sigma,\epsilon}$ and put $\overline{\KK}_R := \KK/R$. Let $Q:V\to\overline{\KK}_R$
be such that
\begin{equation}\label{gquad1}
Q({\bf x}t) ~ = ~ t^\sigma Q({\bf x})t, \hspace{5 mm} (\forall {\bf x}\in V, t\in \KK)
\end{equation}
and there exists a trace-valued $(\sigma,\epsilon)$-sesquilinear form $f_Q:V\times V\rightarrow\KK$ such that
\begin{equation}\label{gquad2}
Q({\bf x}+{\bf y}) ~ = ~ Q({\bf x})+ Q({\bf y}) + (f_Q({\bf x},{\bf y})+ R), \hspace{5 mm} (\forall {\bf x}, {\bf y}\in V).
\end{equation}
We call $Q$ a {\em generalized} $(\sigma,\epsilon)$-{\em quadratic form} ({\em generalized pseudoquadratic form} if we prefer not to mention the pair $(\sigma, \epsilon)$). Of course, for $R = \KK_{\sigma, \epsilon}$ the above definition yields back just pseudoquadratic forms as defined in the previous subsection. If $R \supset \KK_{\sigma, \epsilon}$ then we say that $Q$ is a {\em proper} generalized pseudoquadratic form.

Nearly all we have said for pseudoquadratic forms holds for generalized pseudoquadratic forms as well. In particular, the hypothesis $R\subset \KK$ ensures that the {\em sesquilinearization} $f_Q$ of a generalized pseudoquadratic form $Q$ is uniquely determined by $Q$ (see \cite[Lemma 3.2]{Pas}). Moreover, the existence of $f_Q$ implies that $\overline{R}\subseteq\overline{\KK}^\circ_{\sigma,\epsilon}$, namely $R\subseteq \KK^{\sigma,\epsilon}$ \cite[Theorem 3.3]{Pas}. However, unlike what happens for pseudoquadratic forms, the hypothesis $R \neq \KK$ is not sufficient for $Q$ to be non-trivial. Trivial forms are discussed and characterized in \cite[Proposition 3.5]{Pas}, but we are not going to insist on them here.

Suppose that $Q$ is non-trivial. The radical of $Q$, $Q$-singular vectors of $V$, $Q$-singular points of $\PG(V)$ and totally $Q$-singular subspaces of $V$ or $\PG(V)$ are defined just as for pseudoquadratic forms. A polar space $\Gamma(Q)$ is thus defined, contained in $\Gamma(f_Q)$ as a subgeometry \cite[section 3.3]{Pas}.

As we shall see in the next subsection, alternating forms and generalized pseudoquadratic forms are all we need to describe all embeddings of polar spaces.

\subsubsection{Universal embeddings of embeddable polar spaces}\label{ueps}

Throughout this subsection $\Gamma$ is an embeddable non-degenerate polar space of finite rank $n\geq 2$ and $\varepsilon:\Gamma\to\PG(V)$ is an embedding of $\Gamma$. The underlying division ring of $V$ is denoted by $\KK$ and $Z(\KK)$ is its center.

\begin{theorem}[{\rm Tits \cite[chapter 8]{Tits}}]\label{Tits1}
If $\varepsilon$ is relatively universal then either $\varepsilon(\Gamma) = \Gamma(Q)$ for a non-degenerate pseudoquadratic form $Q$ defined on $V$ or $\mathrm{char}(\KK)\neq 2$ and $\varepsilon(\Gamma) = \Gamma(f)$ for a non-degenerate alternating form $f:V\times V\to\KK$. Moreover, $\varepsilon$ is absolutely universal except in the following two cases:

\begin{itemize}
\item[{\rm (A)}] $n = 2$, $\dim(V) = 4$, $\KK$ is a quaternion division ring and $Q$ is $(\sigma, \epsilon)$-quadratic, where $\sigma$ is the standard involution of $\KK$ and $\KK_{\sigma,\epsilon}$ is a $1$-dimensional subspace of the $Z(\KK)$-vector space $\KK$.
\item[{\rm (B)}] $n = 2$, $\dim(V) = 4$ and $\Gamma$ is a grid of order
  at least $5$ (where by order of $\Gamma$ we mean the cardinality of
  any of its lines). Its $\varepsilon$-image $\varepsilon(\Gamma)$ is a hyperbolic quadric of the projective $3$-space $\PG(V)$.
\end{itemize}
\end{theorem}
No absolutely universal embedding exists for $\Gamma$ in cases (A) and (B). In both cases $\Gamma$ admits more than one embedding, all of them are relatively universal but none of them is homogeneous. We can obtain many of those embeddings from a given one, say $\varepsilon$, by choosing $g\in\mathrm{Aut}(\Gamma)\setminus\mathrm{Aut}_\varepsilon(\Gamma)$ and taking $\varepsilon\cdot g$ as a new embedding (notation as in subsection \ref{lift}). If $\mathrm{Aut}_\varepsilon(\Gamma)g \neq \mathrm{Aut}_\varepsilon(\Gamma)g'$ then $\varepsilon\cdot g \not\cong \varepsilon\cdot g'$. So, the number of isomorphism classes of embeddings of $\Gamma$ is not smaller than the index of $\mathrm{Aut}_\varepsilon(\Gamma)$ in $\mathrm{Aut}(\Gamma)$.

In case (A) the group $\mathrm{Aut}_\varepsilon(\Gamma)$ has index 2 in $\mathrm{Aut}(\Gamma)$. Chosen $\delta\in \mathrm{Aut}(\Gamma)\setminus\mathrm{Aut}_\varepsilon(\Gamma)$, modulo isomorphisms, $\varepsilon$ and $\varepsilon\cdot\delta$ are the unique embeddings of $\Gamma$ (Tits \cite[8.6]{Tits}).

Things go even worse in case (B), not only because the index $[\mathrm{Aut}(\Gamma):\mathrm{Aut}_\varepsilon(\Gamma)]$ increases faster than the order of $\Gamma$ when $\Gamma$ is finite and is infinite when $\Gamma$ is infinite, but also because, when $\Gamma$ is infinite, different fields exist with the same cardinality as the lines of $\Gamma$. Each of these fields yields infinitely many non-isomorphic embeddings of $\Gamma$, which have nothing to do with those associated with another field of the same cardinality.  So, in case (B) with $\Gamma$ infinite, the geometry $\Gamma$ admits no underlying division ring.

Henceforth we assume that $\Gamma$ is not as in cases {\rm (A)} or {\rm (B)} of Theorem {\rm \ref{Tits1}}. So, the embedding $\varepsilon$ is absolutely universal by Theorem {\rm \ref{Tits1}} and, consequently, $\Gamma$ is defined over $\KK$.

\begin{theorem}[{\rm Tits \cite[chapter 8]{Tits}}]\label{Tits2}
Under the previous hypotheses, if $\mathrm{char}(\KK)\neq 2$ then $\varepsilon$ is the unique embedding of $\Gamma$.
\end{theorem}

Moreover, when $\mathrm{char}(\KK) \neq 2$ and $\varepsilon(\Gamma) = \Gamma(Q)$, we can replace the pseudoquadratic form $Q$ with its sesquilinearization $f_Q$.

When $\mathrm{char}(\KK) = 2$, in general $\Gamma$ admits embeddings different from the absolutely universal embedding $\varepsilon$. Explicitly, let $Q$ be a non-degenerate $(\sigma,\epsilon)$-quadratic form such that $\varepsilon(\Gamma) = \Gamma(Q)$ (Theorem \ref{Tits1}). Suppose that $\rad(f_Q) \neq \{0\}$, as it can happen when $\mathrm{char}(\KK) = 2$. For a non-trivial subspace $X$ of $\rad(f_Q)$ we can consider the quotient $\varepsilon_X$ of $\varepsilon$ over $[X]$. Clearly, $\varepsilon_X:\Gamma\to\PG(V/X)$ is an embedding of $\Gamma$ different from $\varepsilon$. Since $\varepsilon$ is absolutely universal, all embeddings of $\Gamma$ different from $\varepsilon$ arise in this way.

With $X$ as above, put $\overline{R}_X := \{Q({\bf x})~|~{\bf x}\in\rad(f_Q)\}$ and $R_X := \{t\in\KK~|~ t+\KK_{\sigma,\epsilon}\in\overline{R}_X\}$ (a subgroup of $\KK$ containing $\KK_{\sigma,\epsilon}$). Suppose firstly that $R_X \neq \KK$ and define $Q_X:V/X\to \overline{\KK}_X := \KK/R_X$ by the clause $Q_X({\bf x}) = Q({\bf x})+R_X$. Then $Q_X$ is a non-degenerate non-trivial generalized $(\sigma,\epsilon)$-quadratic form, its sequilinearization $f_{Q_X}$ is the form induced by $f_Q$ on $V/X\times V/X$ and we have $\varepsilon_X(\Gamma) = \Gamma(Q_X)$ (see \cite{Pas}).

As proved in \cite{Pas}, we have $R_X = \KK$ if and and only if $X = \rad(f_Q)$ and $(\sigma, \epsilon) = (\mathrm{id}_{\KK}, 1)$, namely $\KK$ is a field and $Q$ is quadratic. In this case $f_Q$ induces a non-degenerate alternating form $f_X$ on $V/X\times V/X$ and $\varepsilon_X(\Gamma) = \Gamma(f_X)$.
Note that, in the present situation, if we imitate the previous construction of $Q_X$ then we obtain a useless trivial form.

As it follows from Theorems \ref{Tits1} and \ref{Tits2}, all non-degenerate non-trivial proper generalized pseudoquadratic forms can be obtained as above from a suitable non-degenerate pseudoquadratic form $Q$ in characteristic 2 and a non-trivial subspace $X$ of $\rad(f_Q)$ such that $R_X \subset \KK$.

\begin{note}
When $\rank(\Gamma) > 2$ the existence of the absolutely universal embedding also follows from a criterion of Kasikova and Shult \cite{KS}.
\end{note}

\section{Subspaces spanned by frames}
\label{section3}
\subsection{Definitions and preliminary results}

Let $\Gamma$ be a non-degenerate polar space of finite rank $n \geq 2$. For $2\leq k \leq n$, a {\em partial frame} of $\Gamma$ of {\em rank} $k$ is a pair $\{A,B\}$ of mutually disjoint sets of points, each of size $k$ and such that
\begin{itemize}
\item[(F1)]  $A\subseteq A^\perp$ and $B\subseteq B^\perp$, namely $A$ and $B$ span singular subspaces of $\Gamma$;
\item[(F2)] every point of $A$ is collinear will all but one points of $B$ and, conversely, every point of $B$ is collinear will all but one points of $A$.
\end{itemize}
A {\em complete frame} (also {\em frame} for short) is a partial frame of rank $n$.

\begin{lemma}\label{frame 000}
Let $F = \{A, B\}$ be a partial frame. Then:
\begin{itemize}
\item[{\rm (F3)}] $A$ and $B$ are bases of the projective spaces $\langle A\rangle_\Gamma$ and $\langle B\rangle_\Gamma$ respectively;
\item[{\rm (F4)}] $A^\perp\cap\langle B\rangle_\Gamma = B^\perp\cap\langle A\rangle_\Gamma = \emptyset$.
\end{itemize}
\end{lemma}
\noindent
{\bf Proof.} When $F$ has rank $2$ there is nothing to prove. Suppose that $F$ has rank at least 3 and, by way of contradiction, suppose that $b_1, b_2, b_3\in B$ are collinear. Every point of $A$ is non-collinear with at most one of these points. However they are collinear, hence every point of $A$ is collinear with all of them. Hence $b_i^\perp\supset A$ for $i = 1, 2, 3$, contradicting (F2). Consequently, $B$ is a basis of $\langle B\rangle_\Gamma$. Similarly, $A$ is a basis of $\langle A\rangle_\Gamma$. Claim (F3) is proved.

Turning to (F4), suppose that $A^\perp\cap\langle B\rangle_\Gamma\neq\emptyset$, by way of contradiction. Let $b_0\in A^\perp\cap\langle B\rangle_\Gamma$ and choose a subset $\{b_1, b_2,\dots, b_r\}\subseteq B$ such that $b_0\in \langle b_1, b_2,\dots, b_r\rangle_\Gamma$, minimal among the subsets of $B$ with this property. So $b_0\not\in\langle b_2,\dots, b_r\rangle_\Gamma$. Let $a$ be the unique point of $A$ non-collinear with $b_1$. Then $b_2,\dots, b_r \in a^\perp$ by (F2). However  $b_1\in \langle b_0, b_2,\dots, b_r\rangle_\Gamma$, since $b_0\in\langle b_1, b_2,\dots, b_r\rangle_\Gamma$ but $b_0\not\in\langle b_2,\dots, b_r\rangle_\Gamma$. Therefore $a\perp b_1$, since all points $b_0, b_2,\dots, b_r$ are collinear with $a$. We have reached a contradiction, which forces us to conclude that $A^\perp\cap\langle B\rangle_\Gamma =  \emptyset$. Similarly, $B^\perp\cap\langle A\rangle_\Gamma = \emptyset$. Claim (F4) is also proved. \hfill $\Box$

\bigskip

Every partial frame is contained in complete frame. More precisely:

\begin{lemma}\label{frame 001}
For every partial frame $F = \{A,B\}$ there exists a complete frame $F' = \{A',B'\}$ such that $A' \supseteq A$ and $B' \supseteq B$.
\end{lemma}
\noindent
{\bf Proof.} Let $M$ be a maximal singular subspace containing $A$. Then $M\cap\langle B\rangle_\Gamma = \emptyset$ by (F4). So, we can choose a maximal singular subspace $N \supseteq B$ disjoint from $M$. As $A$ and $B$ are independent sets in the projective spaces $M$ and $N$ (by (F3)), we can extend them to bases $A'$ and $B'$ of $M$ and $N$ respectively, chosen in such a way that (F2) holds for them.  \hfill $\Box$

\bigskip

Let $F = \{A,B\}$ be a partial frame of rank $k$. By (F2), the non-collinearity relation induces a bijection between $A$ and $B$.  Henceforth, when writing $A = \{a_1, a_2,\dots, a_k\}$ and $B = \{b_1, b_2,\dots, b_k\}$ we will always understand that the points $a_1,\dots, a_k$ of $A$ and $b_1,\dots, b_k$ of $B$ are matched in such a way that $a_i\perp b_j$ if and only if $i\neq j$.

By (F3), the singular subspaces $X := \langle A\rangle_\Gamma$ and $Y := \langle B\rangle_\Gamma$ have rank $k$ and (F4) implies that  $X^\perp\cap Y = Y^\perp\cap X = \emptyset$. In particular, $X\cap Y = \emptyset$.

Given a partial frame $F = \{A, B\}$, we put $\langle F\rangle_\Gamma := \langle A\cup B\rangle_\Gamma = \langle X\cup Y\rangle_\Gamma$.

\begin{lemma}\label{frame 00}
Let $F = \{A,B\}$ be a partial frame of rank $k$. Then $\langle F\rangle_\Gamma$ is a non-degenerate polar space of rank $k$.
\end{lemma}
\noindent
{\bf Proof.} Put $\cS := \langle F\rangle_\Gamma$ for short and $X = \langle A\rangle_\Gamma$ and $Y = \langle B\rangle_\Gamma$, as above. Suppose firstly that $k = n$. Then $X$ and $Y$ are mutually disjoint maximal singular subspaces. If $c\in \rad(\cS)$ then $X_c := \langle X, c\rangle_\Gamma$ and $Y_c := \langle Y, c\rangle_\Gamma$ are singular subspaces. However $X$ and $Y$ are maximal. Therefore $X_c = X$ and $Y_c = Y$. It follows that $c\in X\cap Y = \emptyset$; contradiction. Therefore $\cS$ is non-degenerate. Since it contains $X$ and $Y$, which are maximal as singular subspaces of $\Gamma$, necessarily $\rank(\cS) = n$.

Let now $k < n$. By way of contradiction, suppose that $\rad(\cS)\neq\emptyset$ and pick $c\in\rad(\cS)$. Let $M$ and $M'$ be mutually disjoint maximal singular subspaces of $\Gamma$ containing $X$ and $Y$ respectively (they exist since $X^\perp\cap Y = Y^\perp\cap X = \emptyset$). The subspace $Z := M\cap Y^\perp$ is a complement of $X$ in $M$, since $Y^\perp\cap X = \emptyset$ and $\rank(X) = \rank(Y)$. Accordingly, $Z^\perp\supseteq \langle X\cup Y\rangle_\Gamma = \cS$. Hence $c\in Z^\perp$, since $c\in \cS$. Consequently, $c^\perp \supseteq Z^{\perp\perp} = Z$. Similarly, $c^\perp\supseteq Z' := M'\cap X^\perp$. Thus, $c^\perp$ contains both $X$ and $Y$ as well as both $Z$ and $Z'$. Hence $c^\perp\supset M\cup M'$. Therefore $c^\perp\supseteq{\cal R} := \langle M\cup M'\rangle_\Gamma$. However $c\in \cS\subseteq{\cal R}$. It follows that $c\in\rad({\cal R})$. This contradicts what we have already proved in the first part of our proof. So, $\cS$ is non-degenerate.

We shall now prove that $\rank(\cS) = k$. Certainly $\rank(\cS) \geq k$, since $\cS$ contains $X$ and $Y$, which are singular subspaces of rank $k$. By way of contradiction, suppose that $\rank(\cS) > k$. Then $\cS$ admits a
singular subspace $X' \supset X$ of rank $k+1$. As $Y^\perp\cap X = \emptyset$, $X'\cap Y^\perp\neq\emptyset$ is a point, say $c$. Clearly, $c^\perp \supseteq X\cup Y$. Moreover $c\in \cS$. Therefore $c\in \rad(\cS)$. However we have already proved that $\rank(\cS) = \emptyset$. We have reached a final contradiction.  \hfill $\Box$

\begin{lemma}\label{frame 01}
Assume that $\Gamma$ is embeddable and let $\varepsilon:\Gamma\rightarrow\PG(V)$ be an embedding. Suppose that $\Gamma$ is spanned by a frame. Then $\dim(V) = 2n$ and the embedding $\varepsilon$ is relatively universal and admits no proper quotient.
\end{lemma}
\noindent
{\bf Proof.} The geometry $\Gamma$ is spanned by the $2n$ points of a frame. Hence $\dim(V) \leq 2n$. However $\dim(V) \geq 2n$ since $\Gamma$ has rank $n$ and is non-degenerate. Therefore $\dim(V) = 2n$.

Let $\varepsilon':\Gamma\rightarrow\PG(V')$ be such that $\varepsilon'\to\varepsilon$. Then $\dim(V') = 2n$ by the above. As $\dim(V) = \dim(V') = 2n$, necessarily $\varepsilon' \cong \varepsilon$. This shows that $\varepsilon$ is relatively universal. In the same way, starting from $\varepsilon\to\varepsilon'$ instead of $\varepsilon'\to\varepsilon$, we see that $\varepsilon$ admits no proper quotients.   \hfill $\Box$

\begin{corollary}\label{frame 02}
Suppose that $\Gamma$ is embeddable but it is not as in cases {\rm (A)} and {\rm (B)} of Theorem {\rm \ref{Tits1}}. If $\Gamma$ is spanned by a frame, then it admits just one embedding.
\end{corollary}
\noindent
{\bf Proof.} Easy, by combining Lemma \ref{frame 01} with Theorem \ref{Tits1}.  \hfill $\Box$

\subsection{A special case of Theorem \ref{main 1}}

Throughout this section $\Gamma$ is an embeddable non-degenerate polar space of finite rank $n \geq 2$, $\varepsilon:\Gamma\rightarrow\PG(V)$ is a relatively universal embedding, $\KK$ is the underlying division ring of $V$ and $F = \{\{e_1,e_2,\dots, e_k\},\{f_1, f_2,\dots, f_k\}\}$ is a partial frame of $\Gamma$ of rank $k$. We shall prove the following:

\begin{theorem}\label{frame 1}
Put $\langle\varepsilon(F)\rangle_V := \langle \varepsilon(e_1),\dots, \varepsilon(e_k), \varepsilon(f_1),\dots, \varepsilon(f_k)\rangle_V$. Then
\begin{equation}\label{frame 1 eq}
\langle F\rangle_\Gamma ~ = ~ \varepsilon^{-1}(\langle \varepsilon(F)\rangle_V).
\end{equation}
\end{theorem}

In short, the subspace $\langle F\rangle_\Gamma$ arises from the embedding $\varepsilon$. The following lemma is the first step in the proof of Theorem \ref{frame 1}

\begin{lemma}\label{frame 2}
Suppose that {\rm (\ref{frame 1 eq})} holds whenever $k = 2$. Then {\rm (\ref{frame 1 eq})} holds for any $k = 2, 3,\dots, n$.
\end{lemma}
\noindent
{\bf Proof.} To fix ideas, suppose that $\varepsilon(\Gamma) = \Gamma(Q)$ for a pseudoquadratic form $Q$. With this assumption, we miss only the case where $\varepsilon(\Gamma)$ is associated to an alternating form and $\mathrm{char}(\KK) \neq 2$ (see Theorem \ref{Tits1}), but all we are going to say holds for that missing case as well.

Let ${\bf e}_i$ and ${\bf f}_i$ be representative vectors of the points $\varepsilon(e_i)$ and $\varepsilon(f_i)$. By Lemmas \ref{frame 001} and \ref{frame 01}, the set $\{{\bf e}_1,\dots, {\bf e}_k, {\bf f}_1,\dots, {\bf f}_k\}$ is independent in $V$. Put
\[X^+ := \langle e_1,\dots, e_k\rangle_\Gamma, \hspace{4 mm}  X^- := \langle f_1,\dots, f_k\rangle_\Gamma, \hspace{4 mm} S^+ := \langle {\bf e}_1,\dots, {\bf e}_k\rangle_V, \hspace{4 mm} S^- := \langle {\bf f}_1,\dots, {\bf f}_k\rangle_V, \]
\[G ~:=~ \varepsilon(F) ~= ~ \{\{[{\bf e}_1],\dots, [{\bf e}_k]\},\{[{\bf f}_1],\dots, [{\bf f}_k]\}\}.\]
Then $S^\pm$ is totally $Q$-singular, $[S^\pm] = \varepsilon(X^\pm)$, $G$ is a partial frame of rank $k$ in $\Gamma(Q)$ and $\langle\varepsilon(F)\rangle_V = [S^++S^-]$. Property (\ref{frame 1 eq}) amounts to the following:
\begin{equation}\label{frame 2 eq}
\mbox{if } ~ {\bf 0} \neq {\bf v} \in S^++S^- \mbox{ and } ~ Q({\bf v}) = \bar{0} ~ \mbox{ then } ~  [{\bf v}] \in \varepsilon(\langle F\rangle_\Gamma).
\end{equation}
Let ${\bf v} = {\bf v}^+ + {\bf v}^-$ be $Q$-singular with ${\bf v}^\pm \in S^\pm$ and ${\bf v}\neq {\bf 0}$. The vector ${\bf v}^\pm$ is $Q$-singular, since $S^\pm$ is totally $Q$-singular. Moreover, if ${\bf v}^\pm \neq {\bf 0}$ then $[{\bf v}^\pm] = \varepsilon(x^\pm)$ for a suitable point $x^\pm \in X^\pm$. Clearly,  ${\bf v}^+$ and ${\bf v}^-$ cannot be both null, since ${\bf v}\neq {\bf 0}$ by assumption. Suppose that one of ${\bf v}^+$ or ${\bf v}^-$ is null, say ${\bf v}^- = {\bf 0}$. Then ${\bf v} = {\bf v}^+$, hence $[{\bf v}] = \varepsilon(x^+)$ and we are done.

Assume that none of ${\bf v}^+$ and ${\bf v}^-$ is null. Suppose firstly that $[{\bf v}^+]\perp [{\bf v}^-]$ as points of $\Gamma(Q)$. Then $[{\bf v}^+]$ and $[{\bf v}^-]$ are collinear points of $\Gamma(Q) = \varepsilon(\Gamma)$. Accordingly, $x^+$ and $x^-$ are collinear in $\Gamma$ and $\varepsilon$ maps the line $\ell := \langle x^+, x^-\rangle_\Gamma$ of $\Gamma$ onto the projective line $L := \langle [{\bf v}^+], [{\bf v}^-]\rangle_V$ of $\PG(V)$. As $[{\bf v}]\in L$, necessarily $[{\bf v}] = \varepsilon(x)$ for some point $x\in \ell$. The conclusion of (\ref{frame 2 eq}) holds in this case.

Finally, suppose that ${\bf v}^+\not\perp{\bf v}^-$. We can choose non-zero vectors ${\bf w}^\pm\in S^\pm$ in such a way that $G' := \{\{[{\bf v}^+], [{\bf w}^+]\},\{[{\bf v}^-],[{\bf w}^-]\}\}$ is a frame of rank $2$ in $\Gamma(Q)$. As ${\bf w}^\pm\in S^\pm$, there exist points $y^\pm \in X^\pm$ such that $[{\bf w}^\pm] = \varepsilon(y^\pm)$. We have $x^+\perp y^-$ and $x^-\perp y^+$ because $\varepsilon$ induces and isomorphism from $\Gamma$ to $\Gamma(Q)$. Hence $F' = \{\{x^+,y^+\},\{x^-,y^-\}\}$ is a a frame of $\Gamma$ of rank $2$ and $G' = \varepsilon(F')$. By assumption, (\ref{frame 1 eq}) holds for frames of rank 2. Hence $\varepsilon$ maps $\langle F'\rangle_\Gamma = \langle x^+, y^+, x^-, y^-\rangle_\Gamma$ onto $\langle [{\bf v}^+], [{\bf w}^+], [{\bf v}^-],[{\bf w}^-]\rangle_V\cap \Gamma(Q)$. However $[{\bf v}]\in\langle [{\bf v}^+], [{\bf v}^-]\rangle_V\cap \Gamma(Q)$, since ${\bf v} = {\bf v}^++{\bf v}^-$ and $Q({\bf v}) = \bar{0}$. On the other hand, $\langle F'\rangle_\Gamma \subseteq \langle X^+\cup X^-\rangle_\Gamma = \langle F\rangle_\Gamma$. Therefore $[{\bf v}] = \varepsilon(x)$ for some $x\in\langle F\rangle_\Gamma$, as claimed in (\ref{frame 2 eq}). \hfill $\Box$

\bigskip

In view of Lemma \ref{frame 2}, in order to prove Theorem \ref{frame 1} we only need to prove that (\ref{frame 1 eq}) holds when $k = 2$. Thus, let
$F = \{\{e_1, e_2\},\{ f_1, f_2\}\}$ and let ${\bf e}_1, {\bf e}_2, {\bf f}_1, {\bf f}_2$ be representative vectors of $\varepsilon(e_1), \varepsilon(e_2), \varepsilon(f_1), \varepsilon(f_2)$ respectively. Recall that $[{\bf e}_1]\perp[{\bf e}_2]$, $[{\bf f}_1]\perp[{\bf f}_2]$ and $[{\bf e}_i]\perp [{\bf f}_j]$ (as points of $\varepsilon(\Gamma)$) if and only if $i\neq j$. We put $W := \langle {\bf e}_1, {\bf e}_2, {\bf f}_1, {\bf f}_2\rangle_V$, $\Gamma_0 := \langle F\rangle_\Gamma$ and $\Gamma_1 := \varepsilon^{-1}([W])$. Obviously $\Gamma_0\subseteq \Gamma_1$. We must prove that $\Gamma_0 = \Gamma_1$.

Let $\varepsilon_0:\Gamma_0\to\PG(W)$ and $\varepsilon_1:\Gamma_1\to\PG(W)$ be the restrictions of $\varepsilon$ to $\Gamma_0$ and $\Gamma_1$ respectively, with codomain $\PG(V)$ replaced by $\PG(W) = [W]$. Recall that $\varepsilon(\Gamma)$ is associated with a reflexive sesquilinear form or a pseudoquadratic form, by Theorem \ref{Tits1} and the hypothesis that $\varepsilon$ is relatively universal. Hence the same holds for $\varepsilon_1(\Gamma_1) = [W]\cap\varepsilon(\Gamma)$. Moreover, the form associated to $\varepsilon_1(\Gamma_1)$ is non-degenerate, since $\varepsilon_1(\Gamma_1) = \varepsilon(\Gamma_1)$ is the span in $\varepsilon(\Gamma)$ of the frame $\varepsilon(F)$ of $\varepsilon(\Gamma)$.

The embedding $\varepsilon_0$ is relatively universal by Lemma \ref{frame 01}. Accordingly, by Theorem \ref{Tits1} a non-degenerate reflexive sesquilinear or pseudoquadratic form of $W$ also exists which describes $\varepsilon_0(\Gamma_0)$. Moreover, by Theorem \ref{Tits2}, when $\mathrm{char}(\KK)\neq 2$ we can assume that both $\varepsilon_1(\Gamma_1)$ and $\varepsilon_0(\Gamma_0)$ are associated to reflexive sesquilinear forms. When $\mathrm{char}(\KK) = 2$, pseudoquadratic forms can be chosen for both $\varepsilon_1(\Gamma_1)$ and $\varepsilon_0(\Gamma_0)$ (see Theorem \ref{Tits1}).

\begin{lemma}\label{Frame dispari}
Let $\mathrm{char}(\mathbb{K}) \neq 2$. Then $\Gamma_0 = \Gamma_1$.
\end{lemma}
\noindent
{\bf Proof.} Let $\phi, \psi : W\times W\rightarrow \mathbb{K}$ be reflexive sesquilinear forms such that $\varepsilon_1(\Gamma_1) = \Gamma(\phi)$ and $\varepsilon_0(\Gamma_0) = \Gamma(\psi)$. Let $(\sigma, \epsilon)$ and $(\rho, \eta)$ be the admissible pairs associated with $\phi$ and $\psi$ respectively. We know that if $\phi$ is not an alternating form then we can assume that $\epsilon = 1$. Similarly, we can assume that either $\eta = 1$ or $(\rho, \eta) = (\mathrm{id}_\KK, -1)$. In any case, $\epsilon, \eta  \in  \{1, -1\}$.

Let $\cQ := \langle {\bf e}_1, {\bf e}_2\rangle_V\cup\langle {\bf f}_1, {\bf f}_2\rangle_V\cup\langle {\bf e}_1, {\bf f}_2\rangle_V\cup\langle {\bf f}_1, {\bf e}_2\rangle_V$. Two nonzero vectors of the set $\cQ$ are orthogonal if and only if the corresponding points of $\Gamma$ are collinear. Hence $\phi$ and $\psi$ define the same orthogonality relation on the set $\cQ$. Accordingly, both the following hold for any $t,s\in \KK$:
\begin{equation}\label{dispari 1}
\begin{array}{ccc}
\phi({\bf e}_1+{\bf f}_2t, {\bf e}_2 + {\bf f}_1s) = 0 & \mbox{if and only if} & \psi({\bf e}_1+{\bf e}_2t, {\bf f}_1 + {\bf f}_2s) = 0,\\
\phi({\bf e}_1+{\bf e}_2t, {\bf f}_1 + {\bf f}_2s) = 0 & \mbox{if and only if} & \psi({\bf e}_1+{\bf e}_2t, {\bf f}_1 + {\bf f}_2s) = 0.
\end{array}
\end{equation}
Assume to have chosen ${\bf e}_1, {\bf e}_2, {\bf f}_1, {\bf f}_2$ in such a way that $\phi({\bf e}_1, {\bf f}_1) = \phi({\bf e}_2,{\bf f}_2) = 1$, as we can. For $i = 1, 2$ put $c_i := \psi({\bf e}_i, {\bf f}_i)$. By the first of the two equivalences (\ref{dispari 1}) we see that the condition $s + t^\sigma\epsilon = 0$ is equivalent to $c_1s + t^\rho c_2^\rho\eta = 0$, namely $s = -t^\sigma\epsilon$ if and only if $c_1s = -t^\rho c_2^\rho\eta$, which implies
\begin{equation}\label{dispari 2}
c_1 t^\sigma\epsilon ~ = ~ t^\rho c_2^\rho\eta
\end{equation}
for any $t\in \KK$. With $t = 1$ in (\ref{dispari 2}) we get
\begin{equation}\label{dispari 3}
c_1 ~ = ~ c_2^\rho\eta\epsilon^{-1}.
\end{equation}
We can now substitute $c_1$ with its expression (\ref{dispari 3}) in (\ref{dispari 2}). Recalling that $\epsilon, \eta \in \{1, -1\}$ we obtain
$c_2^\rho t^\sigma =  t^\rho c_2^\rho$ for any $t\in \KK$, namely
\begin{equation}\label{dispari 4}
 t^\rho ~ = ~\kappa  t^\sigma \kappa^{-1}
\end{equation}
where $\kappa := c_2^\rho$. However $c_2^\rho = \kappa c_2^\sigma \kappa^{-1}$ by (\ref{dispari 4}). Hence $c_2^\rho = c_2^\sigma$. We now consider the second equivalence of (\ref{dispari 1}). From it we obtain that $1+t^\sigma s = 0$ if and only if $c_1+ t^\rho c_2 s = 0$, namely
$s = -t^{-\sigma}$ (for $t\neq 0$) if and only if $c_1 = t^\rho c_2 t^{-\sigma}$, which yields  

\begin{equation}\label{dispari 6 bis}
c_1 t^\sigma =  t^\rho c_2 \,\,\forall t\in \KK.
\end{equation}


For $t=1,$ we have $c_2 = c_1 =: c$. Note that $c^\sigma = c^\rho$, since $c_2^\sigma = c_2^\rho$. So, (\ref{dispari 3}) and (\ref{dispari 6 bis}) can be rewritten as follows
\begin{equation}\label{dispari 7}
c^\sigma ~ = ~ c^\rho ~ = ~ c\epsilon\eta^{-1}.
\end{equation}
\begin{equation}\label{dispari 8}
 t^\rho ~ = ~ ct^\sigma c^{-1}.
\end{equation}
By  (\ref{dispari 7}) we obtain that $cc^{-\sigma}\epsilon = c(c\eta^{-1}\epsilon)^{-1}\epsilon = \eta$ (recall that $\epsilon, \eta \in \{1, -1\}$). So, $cc^{-\sigma}\epsilon ~ = ~ \eta$. This identity combined with (\ref{dispari 8}) shows that $\psi$ and $\phi$ are proportional: $\psi = c\cdot \phi$ (compare Subsection \ref{rsf}, final paragraphs). Therefore $\Gamma(\phi) = \Gamma(\psi)$, namely $\Gamma_0 = \Gamma_1$.  \hfill $\Box$

\begin{lemma}\label{Frame pari}
Let $\mathrm{char}(\mathbb{K}) =  2$. Then $\Gamma_0 = \Gamma_1$.
\end{lemma}
\noindent
{\bf Proof.} Now we can assume that $\varepsilon_0(\Gamma_0) = \Gamma(Q)$ and $\varepsilon_1(\Gamma_1) = \Gamma(P)$ for suitable non-degenerate pseudoquadratic forms $Q:W\to\overline{K}_{\sigma,\epsilon}$ and $P:W\to\overline{K}_{\rho,\eta}$. Let $f_Q$ and $f_P$ be their sesquilinearizations. As we know from Subsection \ref{pseudo}, we can assume that $\epsilon = \eta = 1$.

Since $\dim(W) = 4$ and $Q$ and $P$ are non-degenerate with $\Gamma(Q)$ and $\Gamma(P)$ of rank 2, the forms $f_Q$ and $f_P$ are non-degenerate. The proof of Lemma \ref{Frame dispari} can be repeated for them word by word. We obtain that $f_Q$ and $f_P$ are proportional, say $f_P  = c\cdot f_Q$ and $t^\rho = ct^\sigma c^{-1}$ for a suitable $c\in\KK^*$ such that $cc^{-\sigma} = 1$, namely $c = c^\sigma$ (compare the final part of the proof of Lemma \ref{Frame dispari} and recall that now $\epsilon = \eta = 1$ by assumption). Accordingly, $c^{-1}\cdot P$ is a $(\sigma, 1)$-quadratic form, proportional to $P$ and with the same sesquilinearization as $Q$. So, modulo replacing $P$ with $c^{-1}\cdot P$, we can assume to have chosen $P$ in such a way that $\rho = \sigma$ and $f_P = f_Q =: f$, say.

It remains to prove that, under this assumption, we have $P = Q$. As in the proof of Lemma \ref{Frame dispari}, we assume to have chosen the vectors ${\bf e}_1, {\bf e}_2, {\bf f}_1$ and ${\bf f}_2$ in such a way that $f({\bf e}_i,{\bf f}_i) = 1$ for $i = 1, 2$. Hence $f({\bf e}_is, {\bf f}_it) = s^\sigma t$ for every choice of $s, t\in \KK$. Consequently, recalling that all vectors ${\bf e}_1, {\bf e}_2, {\bf f}_1$ and ${\bf f}_2$ are singular for both $Q$ and $P$, we have
\[Q({\bf e}_is+{\bf f}_it) ~ = ~ P({\bf e}_is + {\bf f}_it) ~= ~ s^\sigma t+\KK_{\sigma,1}\]
for every choice of $s, t\in \KK$ and $i = 1, 2$. In other words, $Q$ and $P$ coincide on the (non-singular) $2$-subspaces $L_1 = \langle{\bf e}_1, {\bf f}_1\rangle_V$ and $L_2 = \langle {\bf e}_2, {\bf f}_2\rangle_V$. However $W = L_1+L_2$. Therefore, every vector ${\bf v}\in W$ is the sum ${\bf v} = {\bf v}_1+{\bf v}_2$ of a vector ${\bf v}_1\in L_1$ and ${\bf v}_2\in L_2$. As $Q$ and $P$ coincide on $L_1$ and $L_2$, we obtain
\[\begin{array}{ccccc}
Q({\bf v}) & = &  Q({\bf v}_1) + Q({\bf v}_2) + (f({\bf v}_1,{\bf v}_2)+\KK_{\sigma,1}) & = & \\
 & = & P({\bf v}_1) + P({\bf v}_2) + (f({\bf v}_1,{\bf v}_2)+\KK_{\sigma,1}) & = & P({\bf v}).
\end{array}\]
Therefore $Q = P$. Hence $\Gamma_0 = \Gamma_1$.   \hfill $\Box$

\bigskip

The proof of Theorem \ref{frame 1} is complete.

\section{Proof of Theorem \ref{main 1}}
\label{section4}
Henceforth $\Gamma$ is an embeddable non-degenerate polar space of rank $n \geq 2$ and $\varepsilon:\Gamma\rightarrow\PG(V)$ is a relatively universal embedding of $\Gamma$. Note that our hypotheses in this section are slightly more general than those of Theorem \ref{main 1}: we do not ask $\varepsilon$ to be absolutely universal.

\subsection{The non-degenerate case}

Let $\cS$ be a proper subspace of $\Gamma$. We assume that, regarded as a polar space, $\cS$ is non-degenerate of rank $m \geq 2$. Of course, $m \leq n$. We shall prove the following:

\begin{theorem}\label{nondeg 1}
The subspace $\cS$ arises from $\varepsilon$.
\end{theorem}
Clearly, $\cS$ contains a partial frame $F = \{A,B\}$ of rank $m$, but in general $\langle F\rangle_\Gamma \subset \cS$. We say that a subset $C$ of $\cS$ is a {\em generating supplement} of $F$ with respect to $\cS$ (a $\cS$-{\em generating supplement} of $F$, for short) if $\langle A\cup B\cup C\rangle_\Gamma = \cS$. Note that we do not assume that $C$ is minimal with respect to this property, not even that $C\cap(A\cup B) = \emptyset$. Having defined $\cS$-generating supplements in this way, they obviously exists: for want of nicer choices, even the full set of points of $\cS$ can be taken as a generating supplement of $F$.

Given a $\cS$-generating supplement $C$ of $F$, let $\omega$ be an ordinal number of cardinality $|\omega| = |C|$ and choose a well ordering $(c_\delta)_{\delta < \omega}$  of $C$. For every ordinal number $\gamma \leq \omega$ define $\cS_\gamma :=  \langle A\cup B\cup \{c_\delta\}_{\delta < \gamma}\rangle_\Gamma$. In particular, $\cS_0 = \langle F\rangle_\Gamma$ and $\cS_\omega = \cS$.

We shall prove by (transfinite) induction that for every $\gamma \leq \omega$ the subspace $\cS_\gamma$ arises from $\varepsilon$. In particular, $\cS = \cS_\omega$ arises from $\varepsilon$, as claimed in Theorem \ref{nondeg 1}.

Theorem \ref{frame 1} provides the initial step of the induction. Indeed $\cS_0 = \langle F\rangle_\Gamma$. The next lemma provides the inductive step from $\gamma$ to $\gamma+1$.

\begin{lemma}\label{nondeg 2}
Let $\cX$ be a subspace of $\Gamma$ such that $\cX$, regarded as a polar space, is non-degenerate of  rank $m \geq 2$ and $\cX$ arises from $\varepsilon$. Let $a$ be a point of $\Gamma$ not in $\cX$ such that $\cY := \langle \cX, a\rangle_\Gamma$ still has rank $m$. Then $\cY$ is non-degenerate and arises from $\varepsilon$.
\end{lemma}
\noindent
{\bf Proof.} Given a point $c\in \rad(\cY)$ and a maximal subspace $M$ of $\cX\subset \cY$, the point $c$ cannot belong to $M$, since $\cX$ is non-degenerate. Hence $\langle M, c\rangle_\Gamma$ is a singular subspace of $\cY$ of rank $m+1$; a contradiction with the hypothesis $\rank(\cY) = m$. So, $\cY$ is indeed non-degenerate.

Put $\cZ := \varepsilon^{-1}(\langle\varepsilon(\cY)\rangle_V)$. We firstly prove that $\rank(\cZ) = m$ and, consequently, $\cZ$ is non-degenerate by the same argument as in the previous paragraph.

By way of contradiction, suppose that $\rank(\cZ) > m$. Let $M$ be a maximal singular subspace of $\cX$ and let $M' \supset M$ be a singular subspace of $\cZ$ properly containing $M$. If $a^\perp\supseteq M$ then $\langle M, a\rangle_\Gamma$ is a singular subspace of $\cY$ of rank $m+1$. This contradicts the hypothesis that $\rank(\cY) = m$. Therefore $a^\perp\not\supseteq M$. Accordingly, $a^\perp\cap M'$ is a hyperplane in the projective space $M' \supset M$ and $M\not\subseteq a^\perp\cap M'$. In particular, $a\not\in M'$.

Choose a point $b\in a^\perp\cap M'\setminus M$. So, $\ell := \langle a, b\rangle_\Gamma$ is a line of $\Gamma$. As $\cX$ arises from $\varepsilon$ and $a\not \in \cX$, the subspace $\cW := \langle\varepsilon(\cX)\rangle_V$ is a hyperplane of $\cW' := \langle \varepsilon(\cX), \varepsilon(a)\rangle_V = \langle\varepsilon(\cY)\rangle_V$.  Accordingly, the projective line $\varepsilon(\ell)$, which is contained in $\cW'$, meets $\cW$ in a point $\varepsilon(c)$, for $c\in \ell$. However $\cX = \varepsilon^{-1}(\cW)$, since $\cX$ arises from $\varepsilon$. Therefore $c\in \cX$. Hence $c\neq a$. Consequently $\ell = \langle a, c\rangle_\Gamma$. It follows that $b\in \langle \cX, a\rangle_\Gamma = \cY$. However $b\in M'\setminus M$ and $M'\supset M$. So, $\cY$ contains a singular subspace of rank $m+1$, namely $\langle M, b\rangle_\Gamma$. We have reached again a contradiction with the hypothesis $\rank(\cY) = m$, which forces us to conclude that $\rank(\cZ) = m$.

The equality $\cZ = \cY$ remains to be proved. By way of contradiction, let $\cY \subset \cZ$ and choose a point $b\in\cZ\setminus\cY$. By the same argument as in the previous paragraph we see that $b\not\perp a$. We shall prove that $b^\perp\cap M = a^\perp\cap M$ for every maximal
singular subspace $M$ of $\cX$.

By way of contradiction, suppose that there exists a maximal singular subspace $M$ of $\cX$ such that $a^\perp\cap M \neq b^\perp\cap M$ and put $M' := \langle b, b^\perp\cap M\rangle_V$. Both $M$ and $M'$ are maximal
singular subspaces of $\cZ$, since $\rank(\cZ) = m$. We have $a\not\in M'$, since $a\not\perp b$. Therefore $a^\perp\cap M'$ is a hyperplane of $\cZ$. As we have assumed that $a^\perp\cap M \neq b^\perp\cap M$, the subspace $a^\perp\cap b^\perp\cap M$ has codimension $2$ in $M'$. Accordingly, $a^\perp\cap b^\perp\cap M \subset a^\perp\cap M'$.

Choose a point $c\in (a^\perp\cap M')\setminus M$. Then $\ell := \langle a, c\rangle_\Gamma$ is a line of $\cZ$. With $\cW$ as in the previous paragraph, the projective line $\varepsilon(\ell)$ meets $\cW$ in a point $\varepsilon(d)$ for some point $d\in\ell$. However $\cX = \varepsilon^{-1}(\cW)$. Hence $d\in \cX \subset \cY$. Moreover, $d\neq a$ since $a\not\in\cX$. Therefore $\ell = \langle a, d\rangle_\Gamma \subset \cY$. Accordingly, $c\in \cY$. However $c\in M' \setminus M$ and $M'\cap M$ is a hyperplane of $M'$. Therefore $M' = \langle M\cap M', c\rangle_\Gamma$. It follows that $M' \subset \cY$. In particular, $b\in \cY$, whereas $b\in\cZ\setminus \cY$ by assumption. This contradiction forces us to conclude that $b^\perp\cap M = a^\perp\cap M$ for every maximal singular subspace $M$ of $\cX$, as claimed. As $b$ is an arbitrary point of $\cZ\setminus \cY$, this property holds for any point $x\in \cZ\setminus \cY$.

With $b$, $M$ and $M'$ as above, let $x$ be an arbitrary point of $M'\setminus(b^\perp\cap M)$, different from $b$. As $b^\perp\cap M$ is a hyperplane of $M'$, the line $\langle b, x\rangle_\Gamma$ meets $b^\perp\cap M$ in a point. Therefore $x\not\in \cY$, since $M \subset \cY$ but $b\not\in \cY$. Accordingly, $M'\cap\cY = M'\cap M = b^\perp\cap M = a^\perp\cap M$ and $x^\perp\cap N = a^\perp\cap N = b^\perp\cap N$ for every point $x\in M'\setminus M$ and every maximal singular subspace $N$ of $\cX$, according to the conclusions of the previous paragraph. We can choose $N$ such that $N\cap M = \emptyset$. Then $N\cap M' = \emptyset$ since $N\cap M'\subseteq M'\cap\cY = M'\cap M$.

With $N$ chosen in this way and $X := a^\perp\cap N$ ($ = x^\perp\cap N$ for any $x\in M'\setminus M$) we obtain that $X^\perp\cap M'$ is equal to $M'\setminus M$, which of course is not a singleton. However $X$ is a hyperplane of $N$ and $N\cap M' = \emptyset$. Hence $X^\perp\cap M'$ is a singleton. We have reached a final contradiction, which compels us to admit that $\cZ = \cY$.  \hfill $\Box$

\paragraph{End of the proof of Theorem \ref{nondeg 1}.} Let $\cS_\gamma$ be defined as at the beginning of this section. Then $\cS_\gamma$ is a subspace of $\cS$, which is non-degenerate of rank $m$, and it contains a given frame of $\cS$. Hence $\rank(\cS_\gamma) = m$. The following remains to be proved, for every $\gamma \leq \omega$.

\begin{itemize}
\item[$(*)$] the subspace $\cS_\gamma$ is a non-degenerate polar space and it arises from $\varepsilon$.
\end{itemize}
We know from Lemma \ref{frame 00} and Theorem \ref{frame 1} that $(*)$ holds true for $\cS_0$. Lemma \ref{nondeg 2} shows that, if property $(*)$ holds for $\cS_\gamma$, and $\gamma < \omega$, then it also holds for $\cS_{\gamma+1}$. It might happen that $c_\gamma\in\cS_\gamma$; if this is the case then $\cS_{\gamma+1} = \cS_\gamma$ and there is nothing to prove.

The case where $\gamma$ is a limit ordinal remains to consider. In this case $\cS_\gamma = \cup_{\delta < \gamma}\cS_\delta$ and $\cS_\delta$ satisfies $(*)$ for every $\delta < \gamma$, by the inductive hypothesis. Let $c\in\rad(\cS_\gamma)$. Then $c\in \cS_\delta$ for some $\delta < \gamma$. Hence $c\in \rad(\cS_\delta)$. However $\cS_\delta$ is non-degenerate by the inductive hypothesis. We conclude that $\rad(\cS_\gamma) = \emptyset$. Turning to the second part of $(*)$, we have $\langle \varepsilon(\cS_\gamma)\rangle_V = \cup_{\delta < \gamma}\langle \varepsilon(\cS_\delta)\rangle_V$. If $\varepsilon(a)\in\langle \varepsilon(\cS_\gamma)\rangle_V$ then $\varepsilon(a)\in\langle \varepsilon(\cS_\delta)\rangle_V$ for some $\delta < \gamma$. As $\cS_\delta$ is supposed to arise from $\varepsilon$, necessarily $a\in \cS_\delta$. Hence $a\in\cS_\gamma$. We have proved that $\cS_\gamma$ also arises from $\varepsilon$.  \hfill $\Box$

\begin{note}
By the same argument as in the proof of Corollary \ref{main 6}, one can see that every $\cS$-generating supplement of $F$ contains a minimal one, necessarily disjoint from $A\cup B$. We could not take this fact into account in our definition of $\cS$-generating supplements because we don't  know how to prove it without using Theorem \ref{nondeg 1}, which at that stage had still to be proved.
\end{note}

\subsection{The degenerate case}

In this subsection $\cS$ is a degenerate subspace of $\Gamma$ such that $\rank_{\mathrm{nd}}(\cS) \geq 2$. The next theorem finishes the proof of Theorem \ref{main 1}.

\begin{theorem}
The subspace $\cS$ arises from $\varepsilon$.
\end{theorem}
\noindent
{\bf Proof.} By Theorem \ref{Tits1}, either $\varepsilon(\Gamma) = \Gamma(Q)$ for a non-degenerate pseudoquadratic form $Q$ of $V$  or $\varepsilon(\Gamma) = \Gamma(f)$ for a non-degenerate alternating form $f:V\times V\to \KK$ with $\mathrm{char}(\KK) \neq 2$. The second case can be dealt with just in the same way as the first one. To fix ideas, we assume that the first case occurs. So, $\varepsilon(\Gamma) = \Gamma(Q)$.

Let $f_Q$ be the sesquilinearization of $Q$. In order to avoid any confusion with collinearity in $\Gamma$, we denote the orthogonality relation defined by $f_Q$ on $V$ by the symbol $\perp_Q$. We use the same symbol for the corresponding orthogonality relation in $\PG(V)$. For instance, for a subset $X$ of $V$, we put $[X]^{\perp_Q} := [X^{\perp_Q}]$.

Let $R := \rad(\cS)$ and put $\cW := \langle \varepsilon(R)^{\perp_Q}\rangle_V$ (a projective subspace of $\PG(V)$). Then $\Gamma(Q)\cap \cW$ spans $\cW$ and $\varepsilon$ induces on $R^\perp$ a projective embedding $\varepsilon_{(R)}:R^\perp\to \cW$.

Let $\Gamma_R$ be the star of $R$ in $\Gamma$, namely the polar space the points and lines of which are the singular subspaces of $\Gamma$ containing $R$ and with rank equal to $r+1$ and $r+2$ respectively, where $r := \rank(R)$. As $\rank_{\mathrm{nd}}(\cS) \geq 2$, we have $n-r \geq 2$. So, $\Gamma_R$ is indeed a (non-degenerate) polar space of rank $n-r \geq 2$. The embedding $\varepsilon_{(R)}$ induces and embedding $\varepsilon_R$ of $\Gamma_R$ in the quotient $\cW_R := \cW/\varepsilon(R)$ of $\cW$ over its projective subspace $\varepsilon(R)$. The $\varepsilon_R$-image $\varepsilon_R(\Gamma_R)$ of $\Gamma_R$ is the star $\Gamma(Q)_{\varepsilon(R)}$ of $\varepsilon(R)$ in $\Gamma(Q)$. A non-degenerate pseudoquadratic form $Q_R$ can be defined on $\cW_R$ in such a way that $\Gamma(Q_R) = \Gamma(Q)_{\varepsilon(R)}$. Accordingly, $\varepsilon_R$ is relatively universal. Indeed $Q_R$ is pseudoquadratic, whereas every non relatively universal embedding of a polar space arises from either a proper generalized pseudoquadratic form or an alternating form in characteristic $2$ (see Subsection \ref{gp} and the final part of Subsection \ref{ueps}).

Let $\cS_R$ be the star of $R$ in $\cS$. Then $\cS_R$ is a non-degenerate subspace of $\Gamma_R$ and $\rank(\cS_R) = \rank_{\mathrm{nd}}(S) \geq 2$. By Theorem \ref{nondeg 1}, the subspace $\cS_R$ arises from $\varepsilon_R$. Consequently, $\cS$ arises from $\varepsilon_{(R)}$, hence from $\varepsilon$ too.  \hfill $\Box$

\section{Subspaces and $3$-dimensional embeddings}
\label{section5}
The (projective) {\em dimension} of an embedding $\varepsilon:\Gamma\to \PG(V)$ of a connected point-line geometry $\Gamma$ is the dimension of $\PG(V)$. Note that, since embeddings are injective and map lines onto projective lines, if a geometry admits a 2-dimensional embedding then any two of its lines have a point in common. Such a geometry is either a projective plane or a pencil of lines (possibly a single line).  Accordingly, a generalized quadrangle admits no 2-dimensional embedding.

\begin{prop}\label{ultima}
For a generalized quadrangle $\Gamma$, suppose that $\Gamma$ admits a $3$-dimensional relatively universal embedding. Then all proper subspaces of $\Gamma$ have non-degenerate rank at most $1$.
\end{prop}
\noindent
{\bf Proof.} Let $\varepsilon:\Gamma\to\PG(V)$ be a relatively universal embedding of $\Gamma$ such that $\dim(\PG(V)) = 3$. As $\rank(\Gamma) = 2$, all subspaces of $\Gamma$ have rank at most $2$ and, if they have non-degenerate rank 2, then they are non-degenerate, namely subquadrangles of $\Gamma$. Let $\cS$ be a subquadrangle of $\Gamma$. The span $\langle\varepsilon(\cS)\rangle_V$ of $\varepsilon(\cS)$ in $\PG(V)$ cannot be a plane or a line, since $\cS$ admits pairs of non-concurrent lines. Therefore $\varepsilon(\cS)$ spans $\PG(V)$. However $\cS$ arises from $\varepsilon$, by Theorem \ref{nondeg 1}.
It follows that $\cS = \Gamma$.  \hfill $\Box$

\bigskip

The quadrangles considered in cases (A) and (B) of Theorem \ref{Tits1} satisfy the hypotheses of Proposition \ref{ultima}. Hence all of their proper subspaces are either pencils of lines or sets of mutually non-collinear points.


\begin{thebibliography}{999}



\bibitem{BC}
F. Buekenhout and A.M. Cohen.
\newblock {\em Diagram Geometry,}
\newblock Springer, Berlin, 2013.


\bibitem{gr} I. Cardinali, L. Giuzzi, A. Pasini,
  \emph{The generating rank of a polar Grassmannian}, to appear in Adv. Geom.

\bibitem{CS90}
 A.M. Cohen and E.E. Shult.
  \newblock {\em Affine polar spaces}
  \newblock Geo. Dedicata {\bfseries 35} (1990), 43-76.


\bibitem{FF}
C.-A. Faure and A. Fr\"{o}licher.
\newblock {\em Modern Projective Geometry,}
\newblock Kluwer, Dordrecht, 2000.

\bibitem{KS}
 A. Kasikova and E.E. Shult.
  \newblock {\em Absolute embeddings of point-line geometries}
  \newblock J. Algebra {\bfseries 238} (2001), 265-291.

\bibitem{Pas90} A. Pasini,
  \newblock {\em On locally polar geometries whose planes are affine},
  \newblock Geo. Dedicata {\bf 34} (1990), 35-36.

\bibitem{Pas} A. Pasini,
  \newblock {\em Embedded polar spaces revisited},
  \newblock IIG, {\bf 15} (2017), 31-72.

\bibitem{PVM} A. Pasini and H. Van Maldeghem
  \newblock {\em Some constructions and embeddings of the Tilde geometry},
  \newblock Note di Matematica, {\bf 21} (2002-2003), 1-33.

\bibitem{Ron} M.A. Ronan
  \newblock {\em Embedding and hyperplanes of discrete geometries},
  \newblock European J. Combin. {\bf 8} (1987), 179-185.

\bibitem{S}
E.~E. Shult.
\newblock{\em Points and Lines,}
\newblock Springer, Berlin, 2010.

\bibitem{Tits}
 J.~Tits.
  \newblock {\em Buildings of Spherical Type and Finite BN-Pairs}
  \newblock Springer Lecture Notes in Math. {\bfseries 386} (1974).
\end{thebibliography}
\end{document}